\documentclass[12pt,a4paper]{amsart}
\usepackage{amsmath,amssymb,a4,color}
\usepackage{pgf,tikz, cases, IEEEtrantools}
\usetikzlibrary{arrows}
\usetikzlibrary{svg.path}
\usepackage{hyperref, enumerate}
\usepackage[left=2cm,right=2cm,top=2cm,bottom=2cm]{geometry}
\usepackage{pgffor}
\usepackage{comment}

\theoremstyle{plain}
\newtheorem{theorem}{Theorem}[section]

\newtheorem{lemma}[theorem]{Lemma}
\newtheorem{corollary}[theorem]{Corollary}

\theoremstyle{definition}

\newtheorem{remark}[theorem]{Remark}

\newtheorem{definition}[theorem]{Definition}

\def\eqn#1$$#2$${\begin{equation}\label#1#2\end{equation}}

\numberwithin{equation}{section}

\title{Regularity of the inverse mapping in Banach function spaces}
\author{Anastasia Molchanova}
\address{Institute for Analysis and Scientific Computing, TU Wien, Wiedner Hauptstraße 8-10, Vienna, Austria and Institute of Mathematics, Acad. Koptyug avenue 4, Novosibirsk, Russia}
\email{anastasia.molchanova@tuwien.ac.at}
\author{Tom\' a\v s Roskovec}
\address{Faculty of Economics, University of South Bohemia, Studentsk\' a 13, \v Cesk\' e Bud\v ejovice, Czech Republic and Faculty of Information Technology, Czech Technical University in Prague, Th\' akurova 9, 160 00 Prague 6, Czech Republic}
\email{troskovec@ef.jcu.cz}
\author{Filip Soudsk\'y}
\address{Faculty of Economics, University of South Bohemia, Studentsk\' a 13, \v Cesk\' e Bud\v ejovice, Czech Republic}
\email{filip.soudsky@tul.cz }

\thanks{The first author was supported by Austrian Science Fund (FWF) project M 2670, the second author was supported by the grant GA\v CR 18-00960Y, and the third author was supported by EF--IGS2017--Soudsk\' y--IGS07P1.}


\begin{document}

\begin{abstract}
We study the regularity properties of the inverse of a bilipschitz mapping $f$ belonging to $W^m X_{\text{loc}}$, where $X$ is an arbitrary Banach function space. Namely, we prove that the inverse mapping $f^{-1}$ is also in  $W^m X_{\text{loc}}$. Furthermore, the paper shows that the class of bilipschitz mappings in $W^m X_{\text{loc}}$ is closed with respect to composition and multiplication.
\end{abstract}
\maketitle

\section{Introduction}
Sufficient conditions, concerning the derivatives, for a $C^k$-smooth mapping in $\mathbb{R}^n$ to be invertible, are provided by the well-known Inverse Function Theorem.
This subject has attracted the attention of many researchers due to a large number of relevant applications. 
There are two main lines of research. 
The first one, motivated by Control Theory, deals with the theorem for mappings in general metric spaces
regarding a variational or alternative formalism, that provides a better fit to practical problems.
For more information on this topic, we refer the interested reader to the
research of Frankowska \cite{Fra1990}, see also \cite{FatSri1998,GayGeoJea2015,Pal1997}, as well as many others not explicitly mentioned here.
The second question appears in connection with PDEs and goes back to Arnold's paper on Hydrodynamics \cite{Arn1966}. 
The technique proposed there rests on an analysis of geodesics belonging to the group of volume-preserving diffeomorphisms of an (orientated) Riemannian manifold. It requires an investigation of the regularity properties other than $C^k$ of the inverse mapping, as well as of the composition of two mappings.
At the same time, concerning Continuum Mechanics, the study of function spaces, different from the ones of smooth or Sobolev mappings, is of great interest. 
In particular, there are advantages in using Sobolev--Orlicz spaces for nonlinear elasticity \cite{Bal1977}, Lorentz spaces for the Shr\"{o}dinger equation \cite{BreGal1980} and for the $p$-Laplace system \cite{AlbCiaSbo2017}, grand Sobolev spaces for $p$-harmonic operators \cite{DonSboSch2013,GreIwaSbo1997}.
Thoroughly studied, has been the question of the regularity of 
derivatives of the inverse mapping. 
Thus, we refer the reader to \cite{HenKos2006} for Sobolev $W^{1,p}$-regularity in the planar case, to \cite{CsoHenMal2010,Hen2008,HenKosOnn2007,Vod2008,Vod2012} for $BV$- and $W^{1,p}$-regularity in spatial case. 
Also, articles \cite{CHO,IncKapTop2013} deal with the regularity of the inverse mapping and the composition of diffeomorphic or bilipschitz $W^{m,p}$-Sobolev mappings.

In this paper, instead of studying the inverse mapping problem for all the classes of function spaces separately, we take a concept that covers all these options at once.
More precisely, we prove a result for the general rearrangement invariant Banach function spaces.
This approach, developed in \cite{BS}, has recently been very fruitful and many authors have considered issues
such as Sobolev embeddings, the regularity of
solutions to given PDEs and so on in this general setting --- see, for example, \cite{AlbCiaSbo2017}.

The inspiration for our research is a result in classical Sobolev spaces from \cite{CHO}, the proof there builds on the classical Sobolev--Gagliardo--Nirenberg inequality. 
This inequality appears in a much more general form in \cite{FioForRosSou2018}, and this allows us to derive the results which follow. 
In the following text, $\beta_X$ stands for the upper Boyd index of a Banach function space $X$ (see Definition~\ref{BoI}).
In what follows, we prove the following three theorems.

\begin{theorem}\label{thm:main}
    Let $m$, $n\in\mathbb{N}$,
    $\Omega$, $\Omega' \subset \mathbb{R}^n$ be open sets, and
    $X$ be a rearrangement invariant Banach function space such that $\beta_X<1$.
    Also, let
    $f \colon \Omega \to \Omega'$ be a locally bilipschitz homeomorphism with
    $f \in W^m X_{\text{loc}} (\Omega, \mathbb{R}^n)$.
    Then it follows that
    $f^{-1} \in W^m X_{\text{loc}} (\Omega', \mathbb{R}^n)$.
\end{theorem}

\begin{theorem}\label{thm:main_composition}
    Let $m$, $n\in\mathbb{N}$, 
    $\Omega$, $\Omega' \subset \mathbb{R}^n$ be open sets, and
    $X$ be a rearrangement invariant Banach function space such that $\beta_X<1$.
    Also, let
    $f \colon \Omega \to \mathbb{R}^n$ be a locally Lipschitz mapping with 
    $f \in W^m X_{\text{loc}} (\Omega, \mathbb{R}^n)$,
    and
    $g \colon \Omega' \to \Omega$ be locally bilipschitz with
    $g \in W^m X_{\text{loc}} (\Omega', \Omega)$.
    Then it follows that
    $f\circ g \in W^m X_{\text{loc}} (\Omega', \mathbb{R}^n)$.
\end{theorem}

\begin{theorem}\label{thm:main_product}
    Let $m$, $n\in\mathbb{N}$, 
    $\Omega \subset \mathbb{R}^n$ be an open set, and
    $X$ be a rearrangement invariant Banach function space such that $\beta_X<1$.
    Also, let
    $f$ and $g \colon \Omega \to \mathbb{R}$ be locally Lipschitz mappings such that
    $f$, $g \in W^m X_{\text{loc}} (\Omega, \mathbb{R})$.
    Then it follows that
    $fg \in W^m X_{\text{loc}} (\Omega, \mathbb{R})$
    and $fg$ is a locally Lipschitz mapping.
\end{theorem}
\begin{remark}
The result for a product of $f$ and $g$ can be even generalized for $f$, $g\colon \Omega \to \mathbb{R}^n$ being mappings and not just functions, then we understand the product $f\cdot g$ as a scalar product and the proof can be done in the same way with the arguments repeated for all coordinates.
\end{remark}

In particular, these theorems are valid for Lorentz and Orlicz spaces.
Since these spaces are of special interest in applications, we provide an explicit formulation for the reader's convenience.

\begin{corollary}\label{cor:Lorentz}
    Let $m$, $n\in\mathbb{N}$,
    $\Omega$, $\Omega' \subset \mathbb{R}^n$ be open sets,
    $p>1$ and $q \geq 1$.
    Also, let
    $$
    \begin{aligned}
        f & \in W^m L^{p,q}_{\text{loc}} (\Omega, \Omega') \text{ be a locally bilipschitz homeomorphism,}\\
        u & \in W^m L^{p,q}_{\text{loc}} (\Omega, \mathbb{R}^n) \text{ be a locally Lipschitz mapping,} \\ 
        \varphi & \in W^m L^{p,q}_{\text{loc}} (\Omega', \Omega) \text{ be a locally bilipschitz mapping, and}\\
        g, h & \in W^m L^{p,q}_{\text{loc}} (\Omega, \mathbb{R}) \text{ be locally Lipschitz mappings.}
    \end{aligned}
    $$
    Then it follows that
    $$
    \begin{aligned}
        f^{-1} & \in W^m L^{p,q}_{\text{loc}} (\Omega', \mathbb{R}^n),\\
        u\circ \varphi & \in W^m L^{p,q}_{\text{loc}} (\Omega', \mathbb{R}^n), \text{ and}\\
        gh & \in W^m L^{p,q}_{\text{loc}} (\Omega, \mathbb{R}).
    \end{aligned}
    $$
\end{corollary}

\begin{remark}
It is well known that in the case $p>1$ and $q \geq 1$ the upper Boyd index $\beta_{L^{p,q}}<1$.
However, for any  $q>1$ the Lorentz space $L^{1,q}$ is not a Banach function space. In fact, it can not be even equivalently renormed. Thus, it needs a different approach, and we leave the case of $L^{1,q}$ open. 
\end{remark}

\begin{corollary}\label{cor:Orlizc}
    Let 
    $\Omega$, $\Omega' \subset \mathbb{R}^n$ be open sets,
    $A$ be a Young function, such that
    there exists a positive constant $c$, for which
    \begin{equation}\label{nepolka}
        \int_0^t\frac{A(s)}{s^2}\,\textup{d}s\leq\frac{A(ct)}{t},
    \end{equation}
    holds for all $t>0$.
    Also, let
    $$
    \begin{aligned}
        f & \in W^m L^{A}_{\text{loc}} (\Omega, \Omega') \text{ be a locally bilipschitz homeomorphism,}\\
        u & \in W^m L^{A}_{\text{loc}} (\Omega, \mathbb{R}^n) \text{ be a locally Lipschitz mapping,} \\ 
        \varphi & \in W^m L^{A}_{\text{loc}} (\Omega', \Omega) \text{ be a locally bilipschitz mapping, and}\\
        g, h & \in W^m L^{A}_{\text{loc}} (\Omega, \mathbb{R}) \text{ be locally Lipschitz mappings.}
    \end{aligned}
    $$
    Then it follows that
    $$
    \begin{aligned}
        f^{-1} & \in W^m L^{A}_{\text{loc}} (\Omega', \mathbb{R}^n),\\
        u\circ \varphi & \in W^m L^{A}_{\text{loc}} (\Omega', \mathbb{R}^n), \text{ and}\\
        gh & \in W^m L^{A}_{\text{loc}} (\Omega, \mathbb{R}).
    \end{aligned}
    $$
\end{corollary}

\begin{remark}
    The inequality \eqref{nepolka} is an equivalent condition to the boundedness of maximal operator and is in fact equivalent to $\beta_{L^{A}}<1$, see \cite{KH} and Remark~\ref{rem:BoI}.
\end{remark}

\section{Preliminaries}

We use the notation $|\cdot|$ for three different operations on three exclusive types of argument. If the argument is of real value, we consider the symbol to be an absolute value. If the argument is matrix or linear operator, we understand the operator norm. If the argument is a set in $\mathbb{R}^n$, we understand $n$-dimensional Lebesgue measure of this set.

In the following text $\Omega$ and $\Omega'$ stand for open subsets of $\mathbb{R}^n$ with finite Lebesgue measure. We denote a \textit{scaling parameter} as
\begin{equation}\label{eta=sp}
\eta:=\frac{|\Omega'|}{|\Omega|}.
\end{equation}

We write
$A(\xi)\lesssim B(\xi)$
if there exists constant $C>0$ independent of  the parameter $\xi$ such that
$
A(\xi)\leq CB(\xi).
$

\subsection{Banach function spaces}

Let us first remind some notions from the theory of Banach function spaces (later in the text referred just as BFS and r.i.\ BFS if the space is also rearrangement invariant).
We refer the reader to \cite{BS} and \cite{FS}
for the theory of BFS.

\begin{definition}
    Given a BFS $X$ and a real number $\alpha>0$,
    the space $X^\alpha$ consists of all measurable mappings $u$ such that 
    $$
    \|u\|_{X^{\alpha}}:=\left(\| |u|^{\alpha}\|_{X}\right)^{\frac{1}{\alpha}} < \infty. 
    $$
\end{definition}
We 
use the 
convention
\begin{equation}\label{ICO}
X^{\infty}=L^\infty.    
\end{equation}

If $\alpha\geq 1$ then $\|\cdot\|_{X^{\alpha}}$ is a Banach function norm (see \cite[\S 1.d]{LT} and \cite{Loz}).
In this case the space $X^\alpha$
is often referred in the literature as an \it $\alpha$-convexification \rm of $X$. 

Consider numbers
$p_i\in[1,\infty]$, $i=1, \dots k$, such that 
$$
\displaystyle{\sum_{i=1}^k}\frac{1}{p_i}=1,
$$
and locally integrable functions $f_i$, $i=1, \dots k$, then
the following \textit{H\" older inequality} 
\begin{equation}\label{eq:HOL}
    \left\|\displaystyle{\prod_{i=1}^k}f_i\right\|_X\leq\displaystyle{\prod_{i=1}^k}\|f_i\|_{X^{p_i}}
\end{equation}
results from \cite[Proposition 1.d.2]{LT} and the induction by $k$.
Let us remind 
the classical \textit{Hardy--Littlewood--Polya principle} \cite[Corollary II.4.7]{BS}.
For an open set $\Omega\subset\mathbb{R}^n$ and a r.i.\ BFS $X(\Omega)$ the following holds. If 
\begin{equation}\label{HLP_classic}
\int_{0}^tf^*(s)\,\textup{d}s\leq \int_0^tg^*(s)\,\textup{d}s\quad\text{ holds for all } 0<t<|\Omega|,
\end{equation}
then 
$$
\|f\|_{X(\Omega)}\leq \|g\|_{X(\Omega)}.
$$

Here $u^*$ is the \textit{non-increasing rearrangement} of a measurable function $u$,
$$
u^*(s):=\inf\{\lambda: |\{|u|>\lambda\}|\leq s\}.
$$

We define also $u^{**}$ for a measurable function $u$ as
$$
u^{**}(s):=\frac{1}{s}\int_0^s u^*(t)\, dt.
$$
The \textit{Luxemburg representation theorem} \cite[Theorem II.4.10]{BS} states that for every r.i.\ BFS $X(\Omega)$ there exists a r.i.\ BFS $\overline{X}(0,|\Omega|)$, referred as a \textit{representation space}, such that
\begin{equation}\label{Lux}
\|f\|_{X(\Omega)}=\|f^*\|_{\overline{X}(0,|\Omega|)}.
\end{equation}
For our purposes we need a more general form of the Hardy--Littlewood--Polya principle, applicable when the underlying measure space is variable. 

\begin{definition}
Let $s$, $a\in(0,\infty)$, 
the \textit{dilation operator} 
$E_s$
is defined on the space of measurable functions on $(0,a)$ by 
$$
E_sf(x):=\begin{cases} f(sx), \quad &\text{for }sx<a, \\
0, \quad &\text{otherwise},
\end{cases}
$$
for all $x>0$.
\end{definition}
Note that for any Banach function space $X$ one has
\begin{equation}\label{NDO}
   \|E_sf\|_X\leq\max\{s^{-1},1\}\|f\|_X. 
\end{equation}
Indeed, it follows from the fact that
$$
\|E_s f\|_{L^1}\leq s^{-1}\|f\|_{L^{1}}, \quad \|E_s f\|_{L^\infty}\leq \|f\|_{L^{\infty}}
$$
and \cite[Theorem 2.2, p.~106]{BS}.

\begin{definition}
Let $\Omega$, $\Omega'\subset\mathbb{R}^n$ be open sets of finite measure, and $X(\Omega)$, $Y(\Omega')$ be a pair of r.i.\ BFS such that
$$
\eta\|f\|_{X(\Omega)}=\|g\|_{Y(\Omega')}
$$
holds providing
$$
E_{\eta}(g^*)=f^*
$$
with respect to notation \eqref{eta=sp}.
Such spaces are called {\it similar spaces}. To unify the notation of all spaces similar to each other, we use the same name for the space independent of the domains, i.e.\ we denote $X(\Omega'):=Y(\Omega')$.
\end{definition}
\begin{lemma}[Hardy--Littlewood--Polya principle for different measure spaces]\label{HLP}
Let $\Omega$, $\Omega'\subset\mathbb{R}^n$ be open sets of finite measure, let $f$ and $g$ be measurable functions on $\Omega$ and $\Omega'$ correspondingly. Let $X(\Omega)$, $X(\Omega')$ be similar r.i.\ Banach function spaces.
If 
$$
\int_0^tE_{\eta}(g^*)(s) \,\textup{d}s\leq \int_0^t f^*(s) \,\textup{d}s\quad \text{holds for all } t\in(0,|\Omega|),
$$
then this implies that
$$
\|g\|_{X(\Omega')}\leq \max\{\eta^{-1},1\}\|f\|_{X(\Omega)}.
$$
\end{lemma}
\begin{proof}
By the Luxemburg representation theorem \eqref{Lux}, estimate \eqref{NDO}
and the classical Hardy--Littlewood--Polya principle \eqref{HLP_classic}
we obtain
$$
\begin{aligned}
\|g\|_{X(\Omega')}&=\|g^*\|_{\overline{X}(0,|\Omega'|)}\\
&\leq\max\{\eta^{-1},1\}\|E_{\eta}g^*\|_{\overline{X}(0,|\Omega|)}\\
&\leq\max\{\eta^{-1},1\}\|f^*\|_{\overline{X}(0,|\Omega|)}\\
&=\max\{\eta^{-1},1\}\|f\|_{X(\Omega)}.
\end{aligned}
$$
\end{proof}

\begin{definition}[Upper Boyd index]\label{BoI}
The \textit{upper Boyd index} of a r.i.\ BFS $X$ is defined by
$$
\beta_X:=\lim\limits_{t\to \infty}\frac{\log(\|E_{1/t}\|_{\overline{X}\to \overline{X}})}{\log t}.
$$
\end{definition}

\begin{remark}\label{rem:BoI}
    Remind that the maximal operator $M$ is bounded on $X$ if and only if the upper Boyd index $\beta_X<1$,
    see \cite[Theorem 1, p.~3]{Sht}, 
    which is a sufficient
    condition for Theorem~\ref{GBFSR_infinity} being valid. The formulas for calculating the Boyd indices of classical function spaces may be found in literature see, for example, \cite{FiKr}.
\end{remark}

\subsection{Some estimates for weak derivatives}

We refer the reader to the classical book \cite{Maz2011} for the theory of Sobolev spaces.
Let $u\colon \Omega \to \mathbb{R}^n$ be a $k$-times weakly differentiable mapping. 
Let us remind, that for almost every fixed $x_0\in\Omega$,
the $k$-th weak derivative $D^ku(x_0)$ is a $k$-linear mapping. It can be represented by a multidimensional matrix or tensor consisting of all weak partial derivatives of $u$ of order $k$.

Let $X(\Omega)$ be a BFS, the \textit{Sobolev space}
$V^kX (\Omega)$
denotes the space of $k$-times weakly differentiable mappings
$u$
with 
$D^k u\in X(\Omega)$. This space is equipped with semi--norm
$$
 \|u\|_{V^k X (\Omega)}:=\left\|D^k u\right\|_{X(\Omega)}<\infty.
$$
The space 
$W^kX(\Omega)$
consists of all $k$-times weakly differentiable mappings
$u$ 
such that 
$$
\|u\|_{W^k X(\Omega)}:=\displaystyle{\sum_{i=0}^k}\|D^i u\|_{X(\Omega)}<\infty.
$$
We also use the notation
$$
W^k X_{\text{loc}}(\Omega):=\{ u\in W^{k}X(G) \text{ for all } G \text{ open and } G\Subset\Omega\},
$$
here and further $G \Subset \Omega$ means that the closure of $G$ is a compact subset of $\Omega$.

\begin{remark}
    For any BFS 
    $X(\Omega)$ 
    one has 
    $X(\Omega)\subset L^1(\Omega)$ 
    provided that
    $|\Omega|<\infty$, 
    which implies 
    $V^k X_{\text{loc}}(\Omega)\subset V^k L^1_{\text{loc}}(\Omega)$, 
    for arbitrary 
    $k\in\mathbb{N}$ 
    and an open $\Omega\subset\mathbb{R}^n$ of finite measure. 
\end{remark}

A mapping
$f\colon \Omega \to \mathbb{R}^n$
is said to be
\textit{locally bilipschitz} if for every ball $B(x_0,\delta)\Subset\Omega$ centered in $x_0$ with radius $\delta$ there exist
$L>0$
such that 
$$
    L^{-1} |x-y| < |f(x) - f(y)| < L |x-y|
$$
holds for all 
$x$, $y \in B(x_0,\delta)$.

\begin{lemma}[{\cite[Corollary 3.19]{AmbFusPal2000}}]\label{lem:uf-1}
    Let 
    $\Omega$, $\Omega' \subset \mathbb{R}^n$ be open
    and 
    $f\in W^{1,1}_{\text{loc}}(\Omega, \mathbb{R}^m)$.
    Suppose that mapping
    $g\colon \Omega' \to \Omega$
    is a bilipschitz homeomorphism.
    Then 
    $f \circ g \in W^{1,1}_{\text{loc}}(\Omega', \mathbb{R}^m)$
    and
    \begin{equation*}
        D\left(f\circ g(y)\right) = Df (g(y)) Dg(y) \quad \text{for almost all } y\in \Omega'.
    \end{equation*}
\end{lemma}

The crucial part of this paper is the Sobolev--Gagliardo--Nirenberg interpolation inequality, which 
enables estimates to be made of lower order derivatives of the function in terms of higher-order ones and the function itself.
Namely, the inequality
$$
    \|D^j u\|_{X}\lesssim \|D^k u\|^{j/k}_{Y}\|u\|^{1-j/k}_{Z},
$$
which was originally stated by Gagliardo \cite{Ga} and Nirenberg \cite{Ni} in case of 
$X$, $Y$, $Z$ being Lebesgue spaces. 
For our purposes, the particular case of the inequality for BFS recently proved in \cite{FFRS} is needed. For the reader's convenience, let us state the theorem here.

\begin{theorem}[Gagliardo--Nirenberg inequality for r.i.\ BFS]\label{GBFSR}
If $j$, $k\in\mathbb{N}$, $1\leq j<k$, and if $X$, $Y$ are rearrangement invariant Banach function spaces over $\mathbb{R}^n$ such that 
$$
Y^{\frac{k}{j}}\stackrel{\textup{loc}}{\hookrightarrow} X,
$$
then the estimate 
\begin{equation}\label{ME}
\|D^j u\|_{X}\lesssim\|(D^k u)^{**}\|_{\raise -4pt \hbox{${}_Y$}}^{\frac{j}{k}}
\|u^{**}\|_{((Y^{\frac{k}{j}})^{X})^{1-\frac{j}{k}}}^{1-\frac{j}{k}}
\end{equation}
holds for all $k$-times weakly differentiable functions $u$ with a constant independent of $u$.
\end{theorem}
As a corollary we obtain the following theorem once we realise that $X^X=L^\infty$.
\begin{theorem}\label{GBFSR_infinity}
    Let $1\leq j<k$ be natural numbers, 
    and $Y$ be a r.i.\ BFS, such that the upper Boyd index $\beta_Y<1$.
    Then the estimate 
    \begin{equation}\label{GNSCa}
        \|D^{j} u\|_{Y^\frac{k}{ j}}\lesssim\|D^{k} u\|_{Y}^{\frac{j}{k}}\|u\|_{L^{\infty}}^{1-\frac{j}{k}}
    \end{equation}
    is valid for all $k$-times weakly differentiable functions $u$.
\end{theorem}
\begin{remark}
    In the following proof we use the notation for BFS $Z=X^Y$, which means that $Z$ is an optimal space such that the H\"older-type inequality $\|fg\|_X\lesssim \|f\|_Y\|g\|_Z$ holds (see \cite[Lemma 2.2]{FFRS}). This tool may be called the space of H\"older multipliers, see \cite{KoLM} for more details.
\end{remark}
\begin{proof}[Proof of Theorem \ref{GBFSR_infinity}]
Let us set $X:=Y^{k/j}$. Note that the assumptions of  Theorem~\ref{GBFSR} are satisfied since $$Y^{k/j}=X\,\Rightarrow\, Y^{k/j}\stackrel{\text{loc}}{\hookrightarrow}X$$ 
holds and thus, from  the H\"older inequality \eqref{eq:HOL}, one has
$$
L^\infty= ((Y^{k/j})^X)^{1-j/k}.
$$
Using Theorem~\ref{GBFSR}
and the convention~\eqref{ICO}, we derive
$$
        \|D^{j} u\|_{Y^\frac{k}{ j}}\lesssim\|(D^{k} u)^{**}\|_{Y}^{\frac{j}{k}}\|u^{**}\|_{L^{\infty}}^{1-\frac{j}{k}}.
$$
The boundedness of the maximal operator on $Y$ is guaranteed by the assumption on the Boyd index in $Y$, it implies $\|(D^ku)^{**}\|_Y\lesssim\|D^ku\|_Y$, the similar property $\beta_{L^\infty}<1$ results in $\|u^{**}\|_{L^\infty}\lesssim \|u\|_{L^\infty}$. By this we deduce \eqref{GNSCa}. 
\end{proof}
\begin{remark}
    In the case of $X=L^\infty$, Theorem \ref{GBFSR_infinity} coincides with the classic case known as the Kolmogorov--Stein inequality. 
\end{remark}
To get the local version of the theorem above we need an extension operator $\mathcal{E}$, for the construction
of which see
\cite[Theorem~5, p.~181]{STEIN}.
Moreover, the boundedness of the extension operator in the case of classical Sobolev spaces $V^{k}L^p$ was proven there. 
The next theorem for the Sobolev space $V^k X$ follows from the general version \cite[Theorem 4.1]{CiRa}.

\begin{theorem}[On the extension operator]\label{The only one extension operator}
Let $B\subset\mathbb{R}^n$ be a ball and $k\in\mathbb{N}$. Then there exists a linear operator, such that for every r.i.\ BFS $X$ it follows that
\begin{enumerate}[(i)]
\item $\mathcal{E}\colon V^k X(B)\rightarrow V^kX(\mathbb{R}^n)$,
\item $\mathcal{E}u|_B=u.$
\end{enumerate}
\end{theorem}

We can now formulate a local Sobolev--Gagliardo--Nirenberg type theorem.

\begin{theorem}\label{GNLOC}
    Let $\Omega\subset\mathbb{R}^n$ be an open set and $1\leq j<k$ be natural numbers. Then for the r.i.\ BFS $X$, with $\beta_X<1$, it follows that
    $$\left( V^k X_{\text{loc}}(\Omega)\cap L^{\infty}_{\text{loc}}(\Omega)\right) \subset V^{j}X^{k/j}_{\text{loc}}(\Omega).$$  
\end{theorem}
\begin{proof}
For $x\in\Omega$ choose a ball $B=B(x,r)\subset\Omega$.
Theorem \ref{The only one extension operator} implies that the extension $\mathcal{E}u$ belongs to $V^kX(\mathbb{R}^n)$. From Theorem~\ref{GBFSR_infinity} we derive
\begin{equation}\label{eq:GNlocal}
\begin{aligned}
\|D^j u\|_{X^{k/j}(B)}&\leq\|D^j( \mathcal{E}u)\|_{X^{k/j}(\mathbb{R}^n)}\\
&\lesssim \|D^k(\mathcal{E}u)\|_{X(\mathbb{R}^n)}^{j/k}\|\mathcal{E}u\|_{L^{\infty}(\mathbb{R}^n)}^{1-j/k}\\
&\lesssim
\|D^k u\|_{X(B)}^{j/k}\|u\|_{L^{\infty}(B)}^{1-j/k}.
\end{aligned}
\end{equation}
The last inequality is valid due to 
the extension operator can be chosen in the way that
$$
    \mathcal{E}\colon L^{\infty}(B)\to L^{\infty}(\mathbb{R}^n)\quad\textup{and}\quad \mathcal{E}\colon V^{k}X(B)\to V^k X(\mathbb{R}^n),
$$
see \cite{CiRa} for details.
\end{proof}

\subsection{High-order derivatives}

We refer the reader to  \cite[\S 10]{ZorII2004}
for the basic properties of multi-linear mappings and differential calculus, which is useful to deal with high-order derivatives.

The critical tool of the paper is the chain rule. 
Formally, for normed vector spaces $S$, $T$, $R$,
mappings $f\colon S \to T$, $g\colon R \to S$ 
and $r \in R$ we compute
\begin{equation*}
    D(f\circ g) (r) \langle h \rangle = Df (g(r))\langle Dg (r) \langle h \rangle\rangle
    \quad \text{for all } h\in R,
\end{equation*}
which can be written in a matrix form as
\begin{equation*}
    D(f\circ g) = (Df \circ g) \cdot Dg.
\end{equation*}
For the second-order derivative we obtain
\begin{equation*}
    D^2(f\circ g) (r) \langle h_1, h_2 \rangle = D^2f (g(r))\langle Dg (r) \langle h_1 \rangle, Dg (r) \langle h_2 \rangle\rangle + Df (g(r)) \langle D^2 g(r) \langle h_1, h_2\rangle \rangle
\end{equation*}
for all 
$h_1$, $h_2\in R$,
which can be expressed in short as 
\begin{equation*}
    D^2(f\circ g) = (D^2f \circ g) \cdot Dg \otimes Dg + (Df \circ g) \cdot D^2g,
\end{equation*}
where $\cdot$ is used to express the composition of (multi-)linear mappings and
$\otimes$ is a tensor product which makes a bilinear mapping from two linear ones, so that composition has sense.
Further,
\begin{equation*}
\begin{aligned}
    D^3(f\circ g) = & (D^3f \circ g) \cdot Dg \otimes Dg \otimes Dg +  (D^2f \circ g) \cdot  D^2g \otimes Dg + \\
    & 2 (D^2f \circ g) \cdot Dg \otimes D^2g  + (Df \circ g) \cdot D^3g.
\end{aligned}
\end{equation*}
Direct calculations show that 
$D^m (f\circ g)$
is made up from terms of the form 
\begin{equation*}
(D^{k_0} f \circ g ) \cdot \bigotimes_{i=1}^{k_0} D^{k_i}g
\end{equation*}
with some coefficients,
where 
$k_o \geq 1$,
$k_i=0$ if and only if $i>k_0$
and
$\sum\limits_{i=1}\limits^{m} k_i= m$.

Moreover,
for multi-linear mappings $A \in \mathcal{L}(S^{l+n},T)$, $B\in \mathcal{L}(R^k,S^l)$, and $C\in \mathcal{L}(R^m,S^n)$ it follows that 
$A\cdot (B \otimes C) \in \mathcal{L}(R^{k+m},T)$ is a $(k+m)$-linear mapping and we can estimate a norm as
$$
|A\cdot (B \otimes C)|\leq |A| |B \otimes C| \leq |A| |B| |C|.
$$

For more details of this topic, we refer the curious reader to \cite{Man2012}, and to \cite{Yok1992} for the tensor calculus.
The corresponding coordinate representation of the high-order chain rule is described in the best possible way in \cite[\S 2.2]{CHO}.
For the sake of simplicity, we will omit 
$\cdot$ and $\otimes$ 
further in the text,
when it can be done without ambiguity.

\section{Proof of Theorem~1.1. The case {$m=2$}}

To start an induction process, we need to investigate the regularity of the second derivative of the inverse mapping.
We start with the Sobolev regularity case.
\begin{theorem}[{Theorem 1.3 of \cite{Hen2008}}]\label{thm:inverse_1_Sobolev}
    Let 
    $\Omega$, $\Omega' \subset \mathbb{R}^n$ be open,
    $p \geq 1$ and suppose that
    $f \colon \Omega \to \Omega'$ is a bilipschitz mapping.
    If $Df \in W^{1,p}_{\text{loc}} (\Omega, \mathbb{R}^{n^2})$,
    then $Df^{-1} \in W^{1,p}_{\text{loc}} (\Omega', \mathbb{R}^{n^2})$.
\end{theorem}

We provide a more general case involving BFS-regularity. 
\begin{theorem}\label{thm:inverse_1_GBFS}
    Let 
    $\Omega$, $\Omega' \subset \mathbb{R}^n$ be open
    and suppose that
    $f \colon \Omega \to \Omega'$ is a bilipschitz homeomorphism.
    Let $X$ be a rearrangement invariant Banach function space. If 
    $Df \in W^{1}X_{\text{loc}} (\Omega, \mathbb{R}^{n^2})$,
    then $Df^{-1} \in W^{1}X_{\text{loc}} (\Omega', \mathbb{R}^{n^2})$.
\end{theorem}

\begin{proof}[Proof of Theorem~\ref{thm:inverse_1_GBFS}]
Since $Df \in W^{1,1}_{\text{loc}} (\Omega, \mathbb{R}^{n^2})$,
by Theorem~\ref{thm:inverse_1_Sobolev} we know that 
$Df^{-1} \in W^{1,1}_{\text{loc}}(\Omega', \mathbb{R}^{n^2})$.
Then, following the proof of \cite[Theorem 1.3]{Hen2008}, we use Lemma~\ref{lem:uf-1} to differentiate the identity
$f \circ f^{-1} = id$
twice to obtain the equation
\begin{equation*}
    D^2f\left(f^{-1}(y)\right) \left(Df^{-1}(y)\right)^2 + Df\left(f^{-1}(y)\right) D^2f^{-1}(y) = 0,
\end{equation*}
Since $f$ is bilipschitz we know also that there exists a positive constant $L$ such that for almost every $y\in\Omega'$ it holds  
\begin{equation*}
    \left|\left(Df\left(f^{-1}(y)\right)\right)^{-1}\right| \leq L, \: |Df^{-1}(y)| \leq L, \: \text{and} \: |J_{f^{-1}}(y)| \geq L^{-n}.
\end{equation*}
and from previous we derive an estimate
\begin{equation}\label{eq:second_D}
\begin{aligned}
    |D^2f\left(f^{-1}(y)\right)| |\left(Df^{-1}(y)\right)^2|  &= |Df\left(f^{-1}(y)\right)| |D^2f^{-1}(y)|\\
        |D^2f^{-1}(y)|   &=|D^2f\left(f^{-1}(y)\right)| |Df\left(f^{-1}(y)\right)|^{-1} |\left(Df^{-1}(y)\right)^2| \\
        &\leq L^3|D^2f\left(f^{-1}(y)\right)|.
\end{aligned}
\end{equation}

Note that $|D^2 f|$ is a measure absolutely continuous with respect to Lebesgue measure (since $|D^2 f| \in L^1_{\text{loc}}$).
Then for chosen $\varepsilon>0$, there exists $\delta>0$ such that $|E|<\delta$ implies
$$
    \int_E |D^2f(x)|\,\textup{d}x<\varepsilon.
$$
    Let $A \subset K \Subset \Omega'$ be measurable and $G \subset K$ be an open set such that $A \subset G$ and
$$
|G\setminus A|<L^{-n}\delta.
$$
By \eqref{eq:second_D} we get, up to multiple of $L$, the following estimate
$$
\begin{aligned}
    \int\limits_{A} |D^2f^{-1}(y)| \, \textup{d}y 
    &\leq \int\limits_{G} |D^2f^{-1}(y)| \, \textup{d}y\\ 
    &\lesssim \int\limits_{G} |D^2f(f^{-1}(y))| |J_{f^{-1}}(y)| \, \textup{d}y.
\end{aligned}
  $$  
By the change-of-variable formula for Lipschitz functions we obtain  
$$    
\begin{aligned}    
\int\limits_{G} |D^2f(f^{-1}(y))| |J_{f^{-1}}(y)| \, \textup{d}y.&= \int\limits_{f^{-1}(G)} |D^2 f (x)| \, \textup{d}x\\
    &=\int\limits_{f^{-1}(A)} |D^2 f (x)| \, \textup{d}x+\int\limits_{f^{-1}(G\setminus A)} |D^2 f (x)| \, \textup{d}x.
\end{aligned}
$$        
Since $f^{-1}(G\setminus A)<L^n L^{-n}\delta$ by the Lipschitz property of $f^{-1}$, the second term can be estimated and, therefore,      
$$
    \int\limits_{A} |D^2f^{-1}(y)| \, \textup{d}y     \lesssim \int\limits_{f^{-1}(A)} |D^2 f (x)| \, \textup{d}x + \varepsilon.
$$    
For the next calculation, set
$   
    \gamma_t:=\min\{L^n t,|f(\Omega)|\}.
$
Recall that 
$$
    \int\limits_{0}\limits^{t} h^{*}(s)\,\textup{d}s 
    = \sup\limits_{ |A|=t} \int\limits_{A} |h(x)| \, \textup{d}x,
$$
where the supremum is taken over all measurable sets $A$ with $|A|=t$.
Then,
$$
\begin{aligned}
    \int\limits_{0}\limits^{t} |D^2 f^{-1}|^{*}(s)\,\textup{d}s 
    & = \sup\limits_{|A| = t} \int\limits_{A} |D^2f^{-1}(y)| \, \textup{d}y  \lesssim \sup\limits_{|A| = t} \int\limits_{f^{-1}(A)} |D^2 f (x)| \, \textup{d}x + \varepsilon
    \\
    &\lesssim \sup\limits_{|A'| = \gamma_t} \int\limits_{A'} |D^2 f (x)| \, \textup{d}x + \varepsilon  = \int\limits_{0}\limits^{\gamma_t} |D^2 f|^{*}(s)\,\textup{d}s + \varepsilon 
    \\
    & \lesssim L^{-n} \int\limits_{0}\limits^{t} E_{\eta}|D^2 f|^{*}(s)\,\textup{d}s + \varepsilon,
\end{aligned}
$$
where $\eta$ is given by \eqref{eta=sp}.

Here, the constant $\varepsilon>0$ can be chosen as small as we wish. Hence, 
\begin{equation*}
    \int\limits_{0}\limits^{t} |D^2 f^{-1}|^{*}(s)\,\textup{d}s 
    \lesssim \int\limits_{0}\limits^{t} E_{\eta}|D^2 f|^{*}(s)\,\textup{d}s,
\end{equation*}
which implies
$$
 \int\limits_{0}\limits^{t} E_{\eta^{-1}}|D^2 f^{-1}|^{*}(s)\,\textup{d}s 
    \lesssim \int\limits_{0}\limits^{t} |D^2 f|^{*}(s)\,\textup{d}s,
$$
holds for all $t>0$.
Then Lemma~\ref{HLP} guarantees that 
$\|D^2 f^{-1}\|_{X(f(\Omega))} \lesssim \|D^2 f\|_{X(\Omega)}$.
\end{proof}


\section{Proof of Theorem~1.1. The case m$\geq 3$\label{sec:m>3}}

The basic idea of the proof follows \cite{CHO} and is to differentiate the identity
$f \circ f^{-1} = id$
to obtain a representation of the second derivative of the inverse mapping
$D^2 f^{-1}$
through the second derivative $D^2 f$ and the first derivatives $Df$ and $Df^{-1}$.
Further, using the Leibniz and chain rules, we represent 
$D^{k} f^{-1}$ as a product of lower order derivatives of $f$ and $f^{-1}$.
Then the Sobolev--Gagliardo--Nirenberg and H\"older inequalities give us a desirable regularity. 

\begin{lemma}[Lemma 3.1 of \cite{CHO}]\label{lem:D2f-1}
    Let 
    $\Omega$, $\Omega' \subset \mathbb{R}^n$ be open.
    Let
    $f \colon \Omega \to \Omega'$ be a bilipschitz homeomorphism
    such that
    $f \in W^{2,1}_{\text{loc}} (\Omega, \mathbb{R}^{n})$.
    Then
    $f^{-1}\in W^{2,1}_{\text{loc}} (\Omega', \mathbb{R}^{n})$
    and 
    \begin{equation}\label{eq:D2f-1}
        D^2 f^{-1}(y) = -Df^{-1}(y) \cdot D^2f(f^{-1}(y)) \cdot (Df^{-1}(y) \otimes Df^{-1}(y))
    \end{equation}
    for almost all 
    $y\in\Omega'$.
\end{lemma}

\begin{lemma} [Lemma 3.3 of \cite{CHO}]
    Let 
    $\Omega$, $\Omega' \subset \mathbb{R}^n$ be open.
    Let
    $f \colon \Omega \to \Omega'$ be a bilipschitz homeomorphism
    such that
    $f \in W^{m,1}_{\text{loc}} (\Omega, \mathbb{R}^{n})$. Then
    $$|D(D^{m-1}f(f^{-1}))|\in L^1_{\text{loc}} (\Omega')$$
    and
    \begin{equation}\label{eq:Dmff-1}
        D(D^{m-1} f (f^{-1}(y))) = D^m f (f^{-1}(y)) \cdot Df^{-1}(y).
    \end{equation}
    for almost all $y\in\Omega'$.
\end{lemma}

\begin{remark}
    Formula~\eqref{eq:Dmff-1} basically means that 
    \begin{equation*}
        D(D^{m-1} f (f^{-1})) = D^m f (f^{-1})\cdot (Df^{-1} \otimes I \otimes \dots \otimes I),
    \end{equation*}
    where $I$ is the identity mapping.
\end{remark}

Since $f^{-1}$ is bilipschitz,
from \cite[Lemma 3.3]{CHO} it is easy to obtain 
\begin{lemma}\label{lem:DD}
    Let 
    $\Omega$, $\Omega' \subset \mathbb{R}^n$ be open. Let $X$ be a r.i.\ BFS with Boyd index $\beta_X<1$.
    Let
    $f \colon \Omega \to \Omega'$ be a locally bilipschitz homeomorphism
    such that
    $f \in W^{m}X_{\text{loc}} (\Omega, \mathbb{R}^{n})$.
    Then
    $$|D(D^{m-1}f (f^{-1}))|\in X_{\text{loc}}(\Omega')$$
    and \eqref{eq:Dmff-1} holds for almost all $y\in\Omega'$.
\end{lemma}

\begin{proof}[Proof of Theorem~\ref{thm:main}]
We will prove the statement using induction on $m$.
The case $m=1$ follows from the fact that $f$ is bilipschitz.
Theorem~\ref{thm:inverse_1_GBFS} ensures the case $m=2$.

Now, consider the general case 
$m\geq 3$.
Assume that 
$|D^k f^{-1}| \in X_{\text{loc}}(\Omega')$ 
results from 
$|D^k f| \in X_{\text{loc}}(\Omega)$
for all 
$1 \leq k \leq m-1$
and any BFS $X$ with $\beta_X<1$.

Again, as in the proof of \cite[Theorem 1.1]{CHO} we differentiate \eqref{eq:D2f-1} $m-2$ times.
We claim that 
$D^m f^{-1}(y)$
is composed of 
\begin{equation}\label{eq:dmf-1}
    D^{k_{-1}}f^{-1}(y) \cdot D^{k_0} f(f^{-1}(y)) \cdot \bigotimes_{i=1}^{k_0} D^{k_i}f^{-1}(y)
\end{equation}
for almost all $y \in \Omega'$.
Here
$k_{-1} \geq 1$,
$k_0 \geq 2$,
$k_i=0$ if and only if $i>k_0$, and
$k_{-1} + \sum\limits_{i=1}\limits^{k_0} k_i= m+1$.

Since $k_{-1}$, $k_i \leq m-1$ for all $i \geq 1$, 
from Theorem~\ref{GNLOC} with 
$u=Df$, $k = m-1$, $j = k_i -1$
we derive that
$$|D^{k_i} f| \in X^{\frac{m-1}{k_i-1}}_{\text{loc}}(\Omega),$$
and hence by the induction assumption we have
$$
|D^{k_i} f^{-1}| \in X^{\frac{m-1}{k_i-1}}_{\text{loc}}(\Omega').
$$
Now, calculate
$$\frac{k_{-1}-1}{m-1} + \frac{k_0-1}{m-1} + \sum_{i=1}^{k_0} \frac{k_i-1}{m-1} = \frac{1}{m-1}\Big(k_{-1}-1 + k_0-1 + \sum_{i=1}^{k_0} (k_i-1)\bigg) =1.$$
Using this equality as indices in inequality \eqref{eq:HOL} 
we have
$$
   \left\|D^m f^{-1}\right\|_X
    \lesssim
   \| D^{k_0} f\|_{X^{\frac{m-1}{k_0-1}}} \|D^{k_{-1}}f^{-1}\|_{X^{\frac{m-1}{k_{-1}-1}}} \prod_{i=1}^{k_0} \|D^{k_i}f^{-1}\|_{X^{\frac{m-1}{k_i-1}}},
$$
which implies \ 
$|D^{m} f^{-1}| \in X_{\text{loc}}(\Omega')$.

\end{proof}


\section{Proof of Theorems~1.2 and 1.3}

We need the next generalization of \cite[Lemma 4.1]{CHO}.
\begin{lemma}\label{lem:DD_composition}
    Let 
    $\Omega$, $\Omega' \subset \mathbb{R}^n$ be open.
    Let
    $g \colon \Omega' \to \Omega$ be a bilipschitz mapping
    with
    $g \in W^{k}X_{\text{loc}} (\Omega', \mathbb{R}^{n})$,
    and 
    $f \in W^{k}X_{\text{loc}} (\Omega, \mathbb{R}^{n})$.
    Then
    $$|D(D^{k-1}f (g))|\in X_{\text{loc}}(\Omega')$$
    and 
    $$
        D(D^{k-1} f (g(y))) = D^k f (g(y)) \cdot Dg(y).
    $$
\end{lemma}
\begin{proof}[Proof of Lemma \ref{lem:DD_composition}]
The proof of the pointwise equality can be carried out in the very same way as in \cite{CHO} since $W^k X_{\text{loc}}(\Omega)\subset W^{k,1}_{\text{loc}}(\Omega)$ and we can use \cite[Lemma 4.1]{CHO}. In order to do so it is enough to realize that $g$ is bilipchitz and thus $Dg$ is bounded. The rest follows from the point-wise equality.
\end{proof}

\begin{proof}[Proof of Theorem~\ref{thm:main_composition}]
Due to the fact that $g$ is bilipschitz Lemma~\ref{lem:uf-1} (applied on $f=u$, $g=F^{-1}$) provides with the case $m=1$.
Lemma~\ref{lem:DD_composition} gives 
$$
    D^2f\circ g = D^2f(g) \cdot (Dg \otimes Dg) + Df(g) \cdot D^2g 
$$
with $|Df(g)|$ and $|Dg|$ bounded a.e.\ and
$|D^2f(g)|$ and $|D^2g|$ belonging to $X_{\text{loc}}(\Omega')$.

Following the proof of \cite[Theorem 1.2]{CHO}, within Lemmata~\ref{lem:DD_composition}, \ref{lem:uf-1} and the Leibniz rule,  we obtain  that $ D^m (f\circ g)(y)$ is composed of
$$
    D^{k_0} f(g(y)) \bigotimes_{i=1}^{k_0} D^{k_i}g(y)
$$
a.e.\ with
$k_0 \geq 1$,
$k_i=0$ if and only if $i>k_0$, and
$\sum\limits_{i=1}\limits^{m} k_i= m$.
Following the same calculations and estimates as for~\eqref{eq:dmf-1} we ensure that 
$D^m(f\circ g) \in X_{\text{loc}}(\Omega')$.
\end{proof}

\begin{proof}[Proof of Theorem~\ref{thm:main_product}]
The Leibniz rule yields  
$$
    D^m(fg)=\sum_{0\leq j \leq m}{\binom{m}{j}}D^j f \otimes D^{m-j} g.
$$
Therefore, it is enough to show that  
$|D^j f \otimes D^{m-j} g|\in X_{\text{loc}}$ 
for all $j$. 
We may exclude the case $j=0$ or $j=m$ 
since these terms are a product of Lipschitz function and function belonging to $X_{\text{loc}}$. 
For now, we exclude the case $j=1$ or $j=m-1$. By the H\"older inequality \eqref{eq:HOL} and the Sobolev--Gagliardo--Nirenberg type estimate \eqref{eq:GNlocal} for both $Df$ and $Dg$ for any ball $B\Subset \Omega$ we obtain 
\begin{equation}\label{eq:4terms}
\begin{aligned}
    \|D^j f \otimes D^{m-j} g\|_{X(B)} &\leq 
    \|D^j f\|_{X^{\frac{m-1}{j-1}}(B)} \|D^{m-j} g\|_{X^{\frac{m-1}{m-j}}(B)} \\
  &\lesssim \|D^{m}f\|^{\frac{j-1}{m-1}}_{X(B)} \|Df\|^{1-{\frac{j-1}{m-1}}}_{L^\infty(B)}\|D^{m-1} g\|^{\frac{m-j}{m-1}}_{X(B)} \|g\|^{1-{\frac{m-j}{m-1}}}_{L^\infty(B)}.
\end{aligned}
\end{equation}
Three out of four terms are finite by assumptions, we estimate the remaining term $\|D^{m-1} g\|^{\frac{m-j}{m-1}}_{X(B)}$ by  \eqref{eq:HOL} and \eqref{eq:GNlocal} as before
$$
\begin{aligned}
    \|D^{m-1} g\|^{\frac{m-j}{m-1}}_{X(B)} &\leq 
    \|D^{m-1} g\|^{\frac{m-j}{m-1}}_{X^{\frac{m-1}{m-2}}(B)} \|1\|_{X^{m-1}(B)} \\
  &\lesssim \|D^{m}g\|^{\frac{m-2}{m-1}}_{X(B)} \|Dg\|^{1-{\frac{m-2}{m-1}}}_{L^\infty(B)}\|1\|_{X^{m-1}(B)}.
\end{aligned}
$$
All four terms of \eqref{eq:4terms} are finite so all the items $D^j f \otimes D^{m-j} g$ belong to the space. In case $j=1$ we estimate the term by \eqref{eq:HOL} and \eqref{eq:GNlocal} as follows
$$
\begin{aligned}
    \|D^1 f \otimes D^{m-1} g\|_{X(B)} &\leq 
    \|D^1 f\|_{X^{\infty}(B)} \|D^{m-1} g\|_{X^{1}(B)} \\
\end{aligned}
$$
The first term can be estimated by $X^{\infty}=L^{\infty}$. The second case can be considered as previous. The case $j=m-1$ is analogous.
\end{proof}
Note that the boundedness of the maximal operator is needed due the application of the Gagliardo--Nirenberg inequality, which was proved so far only in the case of spaces on which the operator is bounded. The case without the boundedness of the maximal operator is still open.




\end{document}